\documentclass{amsart}
\usepackage{amssymb,amsmath,latexsym,enumerate,graphicx}
\usepackage{xcolor}
\usepackage{tikz}

\newtheorem{question}{Question}

\newtheorem{theorem}{Theorem} 
\newtheorem{corollary}[theorem]{Corollary}
\newtheorem{proposition}[theorem]{Proposition}
\newtheorem{lemma}[theorem]{Lemma}

\theoremstyle{definition}
\newtheorem{example}{Example}
\newtheorem{definition}[theorem]{Definition}

\DeclareMathOperator{\cyc}{Cyc}

\def\Z{\mathbf Z}
\def\R{\mathbf R}
\def\N{\mathbf N}
\def\A{{\bf(A)}}
\def\C{{\bf(C)}}
\def\D{{\bf(D)}}
\def\L{{\bf(L)}}
\def\G{{\bf(G)}}
\def\implies{\Longrightarrow}
\def\cat{\cdot}
\def\const{14c}

\title{A central limit theorem for repeating patterns}
\author{Aaron Abrams, Eric Babson, Henry Landau, Zeph Landau, Jamie Pommersheim}
\date{}

\begin{document}

\begin{abstract}
  We prove a central limit theorem for the length of the longest subsequence of a random permutation which follows one of a class of repeating patterns. This class includes every fixed pattern of ups and downs having at least one of each, such as the alternating case considered by Stanley 
  in \cite{stanley-mich} 
  and Widom
  in \cite{widom}. 
  In every case considered the convergence in the limit of long permutations is to normal with mean and variance linear in the length of the permutation.  
\end{abstract}

\maketitle

\section{Introduction}

A famous and celebrated result of Logan-Shepp \cite{logan-shepp} and Vershik-Kerov \cite{vershik-kerov} is that the expected length of the longest increasing subsequence of a random length $n$ permutation is asymptotic to $2\sqrt{n}$. By contrast, Stanley \cite{stanley-mich} showed that the expected length of the longest \emph{alternating} subsequence 
of a random length $n$ permutation is asymptotic to $2n/3$.  The contrast between the square root behavior in the former study and the linear behavior in the latter led Stanley to pose the following question: which expected lengths are possible, for subsequences specified by a repeated pattern of ``ups'' (increases) and ``downs'' (decreases)?  In this paper we show that, with the exceptions of the increasing and decreasing cases, for any prescribed pattern of ups and downs, the longest subsequence of a random length $n$ permutation that follows the given pattern has linear expected length.  Furthermore the distribution of this length obeys a central limit theorem.  Our methods extend to certain more general types of patterns, and we provide some techniques and algorithms for further exploring these types of problems.

To describe the results in more detail, we begin with the well-known ``longest increasing subsequence'' problem.
This problems and its variants have a long history.  For example the Erd\H{o}s-Szekeres theorem published in 1935 \cite{erdos} asserts that any sequence of $n^2+1$ numbers contains either an increasing or a decreasing subsequence of length $n+1$.  This implies that the expected length of the longest increasing subsequence is at least $\frac 1 2 \sqrt{n}$.  Though this growth rate is correct, the constant is not:  in the large $n$ limit the length of the longest increasing subsequence has expectation roughly $2\sqrt{n}$ (as mentioned above, see \cite{logan-shepp,vershik-kerov}) and, by work of Baik-Deift-Johansson \cite{b-d-j}, is distributed according to a Tracy-Widom distribution.  These difficult results are the product of decades of study.

Richard Stanley considered a variant of this problem by looking instead for \emph{alternating} subsequences.  For our purposes a sequence $b_1,\ldots,b_n$ is \emph{alternating} if the elements satisfy
\[b_1<b_2>b_3<b_4>\cdots\] 
Stanley recognized that the problem of finding the longest alternating subsequence of a given sequence can be broken into pieces, allowing one to patch together solutions of smaller instances of the problem to obtain a solution of a given large problem.  As a result this problem is significantly easier than the increasing case.  He showed (\cite{stanley-mich}, see also Section 8 of his ICM address \cite{icm}) that the expected length of the longest alternating subsequence of a sequence of length $n$ is exactly $(4n+1)/6$ if $n>1$, and shortly thereafter the limiting distribution of this length was shown to be Gaussian by Widom \cite{widom} and others. For some history and applications of alternating sequences, see \cite{arlotto}.

In his 2005 paper \cite{stanley-mich}, Stanley posed a question about other ``up/down'' patterns such as up-up-down, which we denote $UUD$.  A sequence $b_i$ follows this pattern if 
\begin{equation}\label{eq:uud}
    b_1<b_2<b_3>b_4<b_5<b_6>b_7<\cdots
\end{equation}
A general up/down pattern consists of a finite string from the alphabet $\{U,D\}$; in this language, the previously studied examples are the increasing $U$ and alternating $UD$.
Specifically Stanley asked for which up/down patterns there are constants $\mu, c$ such that the expected length of a subsequence following the pattern is asymptotic to $\mu n^c$.

We show here that every up/down pattern has this property.  Again, it was shown in \cite{logan-shepp,vershik-kerov} that the patterns $U$ and $D$ have $\mu=2$ and $c=1/2$.  We prove that every ``non-constant'' pattern, i.e. one containing at least one $U$ and at least one $D$, displays linear growth, i.e., $c=1$, as exemplified by the alternating case for which Stanley also proved $\mu=2/3$.
Moreover we prove the stronger fact, generalizing the alternating case analyzed in \cite{widom}, that for any non-constant up/down pattern $w$, the distribution of the length of the longest subsequence following $w$ obeys a central limit theorem:  as $n$ grows, the distribution over $\sigma\in S_n$ of the length of the longest subsequence of $\sigma$ following $w$ tends to a Gaussian.

We do not have general results about the means of these normal distributions, which correspond to the $\mu$'s.  
However we introduce a dynamical system whose analysis in principle can recover this quantity (and in practice can approximate it).

Our main theorem addresses more general patterns than up/down patterns.  A sequence satisfies a particular up/down pattern if each pair of consecutive elements behaves as the pattern dictates:  either increasing or decreasing.  But one could also look at consecutive triples, for example, and specify a sequence of (desired) order patterns, each being a permutation on 3 letters.  For example a sequence satisfying the inequalities in \eqref{eq:uud} also satisfies the pattern $(1\ 2\ 3), (2\ 3\ 1), (3\ 1\ 2)$ if the additional inequalities $b_{3i+1}<b_{3i-1}$ and $b_{3i+2}<b_{3i}$ hold for all $i$.  To verify this one looks at ``windows'' of size 3 along the sequence $b_i$.

We isolate a combinatorial feature of such patterns that is sufficient to deduce a central limit theorem.  For up/down patterns this feature is also necessary.
The property we seek in a pattern $w$ is essentially the following:
if two sequences both follow the pattern $w$ then the concatenation of those two sequences will also follow $w$, after possibly deleting a bounded number of elements where the sequences are joined.  Such patterns we call {\it combinatorial}. (See Definition \ref{def:comb}.)

Every non-constant up/down pattern $w$ is combinatorial.
The increasing pattern $U$, however, is not:  a long increasing sequence followed by a long but lower increasing sequence can not be combined into a longer increasing sequence, even if the deletion of a bounded number of elements is allowed.

The main idea of the proof is to define, for any combinatorial pattern $w$, a positive probability event called a ``patch'' which, when it occurs, we can use to reset our search.
Loosely speaking, the longest subsequence following $w$ that occurs before a patch can be combined with both the patch itself and also the longest subsequence following $w$ that occurs after the patch, resulting in a longer subsequence following $w$. 

The reason this is helpful is that the search for the longest subsequence following $w$ can be broken up at all occurrences of a patch into smaller searches which are independent of each other.  The results of the small searches can then be combined.  The length of the longest subsequence is therefore the sum of a bunch of iid random variables, setting us up for a central limit theorem. 
The endgame involves an application of Anscombe's theorem, because we do not know exactly how many iid random variables we are summing.  The conclusion is that the limiting distribution is Gaussian.

In the cases of interest we describe a dynamic programming inspired approximation scheme which in general will find an approximate solution, i.e.~a subsequence following $w$ whose length is at least $1-\epsilon$ times the length of the longest subsequence following $w$.  This algorithm runs in $O(n^{k/\epsilon})$ time, where $k$ is the length of the pattern.  
In the special case of up/down patterns we improve this to an $O(n)$ algorithm which returns the longest subsequence exactly, rather than an approximation.

We conclude with some questions.

\section{Notation}\label{sec:notation}

While much of the literature on these problems refers to subsequences of a given permutation of the integers $\{1,\ldots,n\}$, in this work we instead consider sequences of real numbers $s=(s_1,\ldots,s_n)\in \R^n$ with each $s_i$ chosen from Lebesgue measure on the interval $[0,1]$.  
The \emph{order} for such a sequence is the bijection
(permutation) $\sigma$ of $\{1,\ldots,n\}$ that ranks the elements in increasing order; specifically, $s_i<s_j$ if and only if $\sigma(i)<\sigma(j)$.  The $s_i$ are almost surely distinct and the induced distribution over permutations is uniform so the two formulations are equivalent.

In general, we use the ``one-line'' notation $(\sigma(1)\ \sigma(2)\ \cdots\ \sigma(n))$ for a permutation $\sigma$ of the set $\{1,\ldots,n\}$.   We will sometimes view a permutation $\sigma$ as the sequence of integers $(\sigma(1), \sigma(2), \dots, \sigma(n))$ (which has order $\sigma$).  

For sequences $s=(s_1, s_2, \ldots, s_n),\ t=(t_1,\ldots,t_m)$ of real numbers, we'll use the following notation and terminology: 
\begin{itemize}
    \item $|s|=n,$ the length of $s$;
    \item we use open and closed interval notation for consecutive subsequences, e.g., $s[a,b]= (s_a, \dots , s_b)$ and $s[a,b)= (s_a, \dots , s_{b-1})$ and so on;
    \item a consecutive subsequence is also called a \emph{window};
    \item the concatenation is $s\cat t=(s_1, \dots, s_n, t_1, \dots, t_m)$, but we often omit the symbol $\cat$;
    \item similarly we use $s^2=s \cat s$, $s^3$, $s^\infty$, etc.
    \item in the context of a given concatenation $st$, and for a given value $c$, the \emph{$c$-juncture} of $s$ and $t$ is the size $2c$ window $s(n-c,n]\cat t[1,c]$ of $st$;
    \item  $\cyc(s)$ is the cyclically shifted sequence $(s_2,s_3,\ldots,s_n,s_1)$.  The operator $\cyc^i$ shifts a sequence by $i$ terms.
\end{itemize}
When the elements of a sequence are not distinct, as in $s^2$, the sequence does not have a well-defined order.  It is helpful to have a procedure for perturbing sequences to avoid repeated elements.  If $t$ is a finite length sequence, possibly with repeated elements, we define a nearby sequence as follows: let $\delta$ be the minimal distance between distinct elements of $t$ and then define a {\it tie-break} of $t$ to be a sequence $\tilde{t}$ where we perturb equal elements of $t$ by different amounts smaller than $\delta/2$ to render them distinct while leaving the original distinct elements of $t$ unchanged.  This construction ensures that $\tilde{t}$ has distinct elements, and the order of $\tilde{t}[i,j]$ is equal to the order of $t[i,j]$ whenever the latter is defined.

We use $\Z_k$ for the additive group of integers modulo $k$.  Subscripts of $W$'s will generally be interpreted as elements of $\Z_k$.

\begin{definition}[{\bf Patterns, length, window size. Nonconstant, simple}] \label{def:pattern}
    An {\it $r$-pattern $w$ of length $k$} consists of $k$ non-empty subsets $W_1, \dots, W_k \subset S_r$ of the permutation group of $r$ elements with the following property:  for each $i\in \Z_k$ and for each $\pi\in W_i$ there exists $\rho\in W_{i+1}$ such that the order of the last $r-1$ elements of $\pi$ equals the order of the first $r-1$ elements of $\rho$.  (Note $W_{k+1}=W_1$.) We refer to $r$ as the \emph{window size}. 
    
    A pattern is \emph{nonconstant} if $W_i \ne W_j$ for some $i,j$.
    A pattern is \emph{simple} if each $|W_i|=1$.
\end{definition}

\begin{definition}[{\bf Following}] \label{def:follows}
    A sequence $s$ \emph{follows} $w$, or is \emph{$w$-following}, if for all $1 \le i \le |s|-r+1$, the order of $s[i, i+r)$ is contained in $W_i$.
\end{definition}

In the special case $r=2$, each element of each $W_i$ is one of the permutations $(1\ 2)$ or $(2\ 1)$, which we denote by $U$ (pronounced ``up") and $D$ (``down'') respectively.  A word of length $k$ in $U$ and $D$ will denote a $2$-pattern of length $k$ with $W_i$ being the singleton set consisting of the $i$th letter of the word.  With this language, an increasing sequence follows the pattern $U$, and likewise an alternating sequence follows $UD$.  (We follow the opposite convention from Stanley, who begins alternating sequences with a ``down.'')

It rarely happens that the concatenation of two $w$-following sequences is also $w$-following.  There are essentially two reasons for this:  a localized boundary effect (where the sequences are joined) could interfere, or the sequences could have a global incompatibility.  We are primarily interested in those patterns $w$ for which the second of these never occurs.  This leads to the following definition.

\begin{definition}[{\bf Combinatorial, merge}] \label{def:comb}
A pattern $w$ is \emph{combinatorial} with \emph{combinatorial constant} $c$ if for every pair of sequences $s,t$ following $w$, there exists a $w$-following sequence $u$ obtained from $st$ by removing at most $c$ elements from the $c$-juncture of $s$ and $t$.  Such a sequence $u$ is called a \emph{merge} of $s$ and $t$.
\end{definition}

We will discuss this definition and give some examples after stating our main theorem in the next section.

\section{Main results}

Our primary focus is on the distribution $L^w_n$ of the length of the longest $w$-following subsequence of a uniformly random sequence chosen from $[0,1]^n$.  We denote the mean of this distribution by $\mu^w_n$.

Our main theorem is that if $w$ is a combinatorial pattern then the distributions $L^w_n$ satisfy a central limit theorem. Let $\Phi(t)=\frac{1}{\sqrt{2\pi}}\int_{-\infty}^{t}e^{-u^2/2}du$ denote the cumulative distribution function for the standard normal distribution.
\begin{theorem} \label{thm:main} 
If $w$ is combinatorial then there exist $\mu\in\R$ and $0<\sigma\in\R$ with 
\begin{equation} \label{e:1}
\lim_{n\rightarrow\infty} \hbox{Prob}\Big{(}\frac{L^w_n-\mu n}{t\sigma\sqrt{n}}<1\Big{)}=\Phi(t) 
\end{equation}for every $t\in \R$.    
\end{theorem}

\begin{question}
    How fast is the convergence in Theorem \ref{e:1}?
\end{question}

The theorem applies in particular to all $2$-patterns except the constant patterns (of any length) consisting of all $U$'s and the constant patterns consisting of all $D$'s; that is to say, the exceptions are precisely the increasing and decreasing cases treated by the Baik-Deift-Johansson Theorem.

To see this, let $w$ be a $2$-pattern of length $k$. 
Note of course that if $w$ \emph{is} a constant pattern consisting of all ups or all downs, then $w$ is not combinatorial.
We will argue that in all other cases, $w$ is combinatorial.

If $w$ is not simple, i.e., if some $W_j=\{U,D\}$, then suppose $s=(s_1,\ldots,s_n)$ and $t=(t_1,\ldots,t_m)$ follow $w$ and let $J$ be maximal such that $J\equiv j \mod k$ and $J \le n$.  Then the sequence
$s_1,s_2,\ldots,s_J,t_{j+1},t_{j+2},\ldots,t_m$ follows $w$ so $w$ is combinatorial.

If $w$ is simple and nonconstant then without loss of generality we assume $W_1=U$ and we set $j$ to be the smallest index with $W_j=D$.  In other words we assume $w$ begins with $j-1$ ups and then a down.
Suppose $s=(s_1,\ldots,s_n)$ and $t=(t_1,\ldots,t_m)$ follow $w$ and let $J$ be maximal such that $J\equiv j \mod k$ and $J \le n$.
If $s_J > t_{j+1}$ then the sequence
$s_1,\ldots,s_J,t_{j+1},\ldots,t_m$
follows $w$.
If $s_J < t_{j+1}$ then the sequence
$s_1,\ldots,s_{J-1},t_j,\ldots,t_m$
follows $w$.
In either case we have removed at most $c=k+j\le 2k$ elements from the appropriate part of the concatenation $st$ to produce a $w$-following sequence, 
showing that $w$ is combinatorial with constant at most $2k$.

With a little more care the combinatorial constant in the above argument can be reduced to $c\le k$.
We know of no combinatorial pattern where is it necessary to take $c>\max\{k,r\}$.

\begin{corollary}
For every nonconstant $2$-pattern $w$ 
the $\{L^w_n \}$ satisfy a central limit theorem. 
\end{corollary}

In particular, if $w$ is a nonconstant $2$-pattern then the expected length of the longest $w$-following subsequence of a random sequence is asymptotically $\mu n$, where the constant $\mu$ depends only on $w$.
This resolves Problem 2 of Stanley \cite{stanley-mich}.

\begin{question}
    Given a nonconstant $2$-pattern, or more generally any combinatorial pattern, can one determine the mean and variance of the corresponding Gaussian?
\end{question}

Again, Stanley has shown that $\mu=2/3$ for the alternating pattern $w=UD$.  We discuss an alternate perspective on this question in Section \ref{sec:UUD}.

Combinatorial patterns with $r>2$ also exist.
It is worth spending a moment to compare the following examples.
The constant $3$-pattern of length $1$ with $W_1=\{(1\ 2\ 3),(2\ 3\ 1),(3\ 1\ 2)\}$ is combinatorial, as is the $3$-pattern of length $4$ with $W_1=W_4=\{(1\ 2\ 3)\}, W_2=\{(2\ 3\ 1)\}, W_3=\{(3\ 1\ 2)\}$.
These are both straightforward exercises.
The $3$-pattern of length $3$ with $W_1,W_2,W_3$ as in the previous example (but with no $W_4$) is not combinatorial; this pattern is called $w_2$ in Section \ref{sec:questions}, where it is discussed further.

The theorem is proved using the notion of a {\it patch}. This is a collection of $O(\max (c,k,r)^2)$ consecutive points of a random sequence having a certain specific form.  The key concept is that patches break the problem up into smaller independent problems.  Specifically, we show that when such a patch occurs, which is a positive probability event, a longest subsequence following a given pattern can be found by combining a longest one before the patch and a longest one after it.  The problem is thus reduced to a sum of independent events which will lead to the central limit behavior.  As the number of patches and therefore the number of iid variables being summed is also random, we apply Anscombe's Theorem to complete the argument.

\section{Patching}

To define our patches we will utilize a long permutation $\tau$ which will describe the ``shape'' of a patch.  We define $\tau$ explicitly in the appendix and establish various technical properties of $\tau$ in Lemma \ref{l:tau}.  Using these properties of $\tau$ we prove Lemma \ref{l:sewing} below which provides the key condition to making our patching argument work.

\subsection{Tracking.}\label{sec:tracking}

We first define the notion of a sequence $t$ \emph{tracking} a permutation $\sigma$. For this purpose we view the permutation $\sigma \in S_m$ as a sequence $(\sigma(1),\ldots,\sigma(m))$, and we say that $t$ \emph{tracks} $\sigma$ if $t$ follows the $m$-pattern of length $m$ given by $W_i=\{\cyc^{i-1}(\sigma)\}$.
Equivalently, the order of $t[i, i+m)$ should agree with the order of $\sigma^\infty[i, i + m )$, for all $1\leq i \leq |t|-m +1$.  Note that there is an important difference between $t$ tracking $\sigma$ and $t$ following the pattern $\{\sigma\}$; for example if $\sigma=(1\ 2)$ then tracking $\sigma$ is equivalent to being alternating.

\begin{lemma} \label{l:sewing}
Given a combinatorial $r$-pattern $w$ of length $k$ with combinatorial constant $c$ there exists an $m\geq \max (4c,4k,r)$ that is a multiple of $k$ and a permutation $\tau \in S_m$ such that the following holds: 

whenever $t, u$ are $w$-following sequences such that:
\begin{itemize}
    \item the last $2m$ elements of $t$ track $\tau$,
    \item both the first and last $2m$ elements of $u$ track $\tau$,
    \item the $m$-juncture of $tu$ tracks $\tau$,
\end{itemize}
then $t u$ follows $w$.
\end{lemma}

We remark that this lemma is straightforward if $|t|$ is a multiple of $k$.  The difficulty arises in verifying the $w$-following condition for windows in $u$, if $|t|$ is not a multiple of $k$.

\begin{proof} 
We define $\tau$ according to Lemma \ref{l:tau}.

Our goal is to verify that 
\begin{equation} \label{e:2} (t \cat u) [i, i+ r) \in W_{i} 
\end{equation}
for all $i, 1\leq i \leq |t  u| - r +1$.
The fact that $t$ follows $w$ implies (\ref{e:2}) for $1 \le i \le |t|-r +1$.  For $|t|-m +1\le i \le |t|+m$, the tracking hypotheses for $u$ and $t  u$ imply that the order of $(t \cat u) [i, i+ r)$ is the same as the order of the window $m$ units to the left $(t \cat u) [i-m, i-m+ r)$ which we already know is in $W_{i-m}= W_{i}$.

Since $m>r$ it remains only to consider $i >|t|+ m$.
Note that the hypothesis about $u$ implies by condition (3) of Lemma \ref{l:tau} that the corresponding path of $u$ in $H$ only passes through loop vertices of $(\sigma',1)$ where $\sigma'$ is the order of $\tau[1,r]$.  Since the path corresponding to $\tau$ in $H$ passes through {\it all} loop vertices of $(\sigma', 1)$ it follows that for each $i$, there exists $j\leq m$ such that $j \equiv i \mod k$ and the order for $u[i,i+r)$ is equal to the order of $u[j, j+r)$.
Thus the order of $(t\cat u)[|t| +i, i+r)$, which is also the order of $u[i, i +r )$, equals the order of $u[j, j+r)$, hence of $(t\cat u)[|t| + j, |t| + j + r)$.  Because $j \leq m$ the latter order has already been shown to be in $W_{|t| + j}= W_{|t|+i}$.  Thus (\ref{e:2}) holds for all $i, 1\leq i \leq |t \cat u| - r +1$. 
 \end{proof}

\subsection{Patches}
A patch for $w$ is one of a specific collection of length $\const m$ sequences that tracks $\tau$.

\begin{definition}[Patch] \label{def:patch}
Given a combinatorial pattern $w$ with combinatorial constant $c$, and a permutation $\tau \in S_m$ given by Lemma \ref{l:sewing}, divide the interval $[0,1]$ into $m$ equally spaced disjoint intervals $I_1=[0, \frac{1}{m}), I_2 =[\frac{1}{m}, \frac{2}{m}),\ldots, I_k=[\frac{m-1}{m}, 1]$.  Further divide each of these intervals into $\const$ equal pieces and define $I_{j_1, j_2} = [\frac{j_1-1}{m} + \frac{j_2-1}{\const m}, \frac{j_1 -1}{m} + \frac{j_2}{\const m})$ for $1\leq j_1 \leq m$, $1 \leq j_2 \leq \const $.

    Let $B$ be the following subcube of $([0,1]^m)^{\const}$: 
    \[ B=\prod_{j=1}^{\const} \left( I_{\tau(1), j} \times I_{\tau(2), j} \times \cdots \times I_{\tau(m), j} \right). \]
We call $B$ a patching set for $w$ and any element of $B$ a \emph{patch} for $w$.  
\end{definition} 

Note the volume of $B$ is $(\const m)^{-\const m}$, and every patch $b\in B$ has the following properties:
\begin{itemize}

    \item $b$ tracks $\tau$,
    \item $b$ follows $w$,
    \item $b_i< b_{i+m}$ for all $i$.
    
\end{itemize}  
In addition, it is a straightforward verification that the patch has the following helpful property:

\begin{lemma} \label{l:shiftinvariant}
    Given a patch $b= (b_1, \dots, b_N)$ then the order of any subsequence $(b_{i_1}, \dots, b_{i_{\ell}})$ is equal to the order of the subsequence  $(b_{i_1 +m}, \dots , b_{i_{\ell} +m})$  (provided $i_{\ell} + m \leq N = 14cm$).
\end{lemma}

\subsection{Using patches to break the problem into pieces}
In Proposition \ref{l:patch} below, we show the key property of a patch: if a sequence $s$ contains a patch, then there exists a maximal subsequence following $w$ that uses the entire window of length $4m$ located in the middle of the patch.  

Starting with a sequence $s$ containing a patch $b$, let $t$ be a maximal subsequence of $s$ that follows $w$.  Decompose $t=X Y Z$ where $Y$ are the elements of $b$ in $t$.  We will show we can replace $Y$ (while preserving the $w$-following property) with $\overline{Y}$ where $\overline{Y}$ has the same length as $Y$ and includes the middle $4m$ elements of the patch $b$.

Notice that using the combinatorial property of $w$ on the sequence $X b Z$ would produce a subsequence that follows $w$ of length lower bounded by $|X| + |Z| + \const m - 2c$ which implies $|Y| \geq \const m -2c$, in other words $Y$ must include all but possibly $2c$ elements of $b$. 

It will be convenient to describe $Y$ by indicating which of the elements it uses in $b$.   We'll call the {\it location profile} of $Y$ the word in $\{0,1\}^{\const m}$ whose $i$th bit indicates whether the $i$th location of $b$ is included in $Y$ (if $1$) or skipped (if $0$).  

\begin{lemma} \label{p:switch} Using the notation above, suppose $XYZ$ follows $w$ and the location profile for $Y$ is of the form $A I ^2 B I^3 C$ where $A$, $B$, $C$ are each binary words of length a multiple of $m$ and  $I$ is the word $1^m$.
Then the subsequence $\overline{Y}$ with location profile $A I ^3 B I^2 C$ has the property that $X\overline{Y} Z$ follows $w$.
\end{lemma}
\begin{proof} Recall that to follow $w$ is to satisfy a local condition that needs to be verified for the order of every window of size $r$ of $X\overline{Y} Z$ (recall also that $m\ge r$).  The fact that $XYZ$ follows $w$ verifies this local condition directly for all windows of length $r$ within $X\overline{Y}Z$ except some of those that occur within $\overline{Y}$, specifically those that lie within the $I^2BI$ portion of the location profile $A I ^3 B I^2 C=AI(I^2BI)IC$ of $\overline{Y}$.  Consider any such window.  By Lemma \ref{l:shiftinvariant}, the window $m$ elements to the left has the same order; but the latter window is one of the windows in $Y$.  The fact that $XYZ$ follows $w$ then verifies the required local condition for these cases. \end{proof}

\begin{proposition} [Patching Proposition] \label{l:patch}
Suppose we are given a combinatorial $r$ pattern $w$ of length $k$ and a sequence $s$ with a patch $b \in B$.  Let $t$ be a maximal $w$-following subsequence of $s$.  Then there exists a maximal $w$-following subsequence of $s$ that coincides with $t$ outside $b$ and that intersects $b$ in a set that includes the entire length $4m$ window in the middle of $b$.
\end{proposition}

\begin{proof}
Let $t$ be a maximal subsequence of $s$ that follows $w$ with the following added property:  as before, decompose $t=XYZ$ where $Y$ are the elements of $b$ in $t$ and require that $Y$ is such that its location profile contains $I^q$ for $q$ as big as possible.  We shall argue that $q \geq 8c  $ and therefore it must contain the middle $4m$ elements of $b$.  First we remark simply that $q \geq 3$ since the maximal length of a location profile with at most $2c$ locations that are $0$'s and does not have a portion of the form $I^3$ is $(6c+2)m < \const m$. Note that $I^3$ cannot appear in the location profile anywhere outside of the $I^q$ section since we could use Proposition \ref{p:switch} to switch one of the copies of $I$ next to $I^q$ thus violating the maximality of $q$ (technically, Proposition \ref{p:switch} only applies to an $I^3$ located to the right of $I^q$ but the same proof works with left and right reversed). Since there are at most $2c$ skipped values of $0$ in the location profile, the restriction that no $I^3$ appears means that the length of the location profile outside of $I^q$ must be no longer than $6cm$ which implies that the length of $I^q$ is at least  $\const m - 6c m= 8cm$.   Since this length is more than $|b|/2 + 2m = (7c +2) m$ we can conclude that the middle length $4m$ window of $b$ is contained in $Y$.  
\end{proof}

Proposition \ref{l:patch} gives a road map for how to break the problem of finding longest $w$-following subsequences into a collection of smaller problems.  Starting with a given sequence $s$ and a patching set $B$ for $w$, we let $b_1 \in B$ be the first occurrence of a patch in $s$, i.e., a length $14cm$ window of $s$ that is in $B$.
Let $b_2$ be the first patch that starts after the completion of $b_1$, and so on.  Define the {\it patch decomposition} of $s$ to be the decomposition of $s$ into the disjoint pieces $s[1, j_1], s(j_1, j_2], \dots , s(j_{\ell -1} , j_{\ell}], s(j_{\ell},|s|]$ where $j_i, j_i +1$ are the locations of the middle two elements of $b_i$.  We can combine the Patching Proposition and Lemma \ref{l:tau} to reduce finding an approximation to the longest $w$-following subsequence to solving a set of reduced sized questions:
\begin{lemma} \label{l:interval}
Given a sequence $s$ of length $n$ and a patching set $B$, let 
\[s[1, j_1], s(j_1, j_2], \dots , s(j_{\ell -1} , j_{\ell}], s(j_{\ell},n] \]
be its patch decomposition.  Define $u_0$ to be the longest subsequence of $s[1, j_1]$ that follows $w$ and includes the last $2m$ elements of $s[1, j_1]$.  Define $u_i$, $i\geq 1$ to be the longest subsequence of $s(j_i , j_{i+1}]$ that follows $w$ and includes both the first and last $2m$ elements of $s(j_i , j_{i+1}]$.  Then the sequence $u= u_0 u_1 \cdots u_{\ell-1}$ follows $w$ and the length of the longest $w$-following subsequence of $s$ lies in the interval $[|u|, |u| + n - j_{\ell}]$. 
\end{lemma}

\begin{proof}   By Lemma \ref{l:sewing} we have that $u= u_0 u_1 \cdots u_{\ell-1}$ follows $w$ and therefore its size is a lower bound for the maximal length subsequence that follows $w$.\footnote{This bound could be incremented by $2m$ by including the first $2m$ elements of the final interval of the patch decomposition.}
Repeated use of the Patching Proposition implies there exists a longest subsequence $t$ of $s$ that follows $w$ and includes the first and last $2m$ elements of all but the first and last element of the patch decomposition (in addition $t$ can be found to include the last $2m$ elements of $s[1,j_1]$ and the first $2m$ elements of $s[j_{\ell}, n] $).  It follows that the size of $t$ within $s[1, j_{\ell}]$ is bounded by $\sum_{i=0}^{\ell -1} |u_i|$ since each $u_i$ was chosen of maximal possible length.  The result follows by noting that  $|t| \leq \sum_{i=0}^{\ell -1} |u_i| + |s(j_{\ell}, n]| =  |u| + n - j_{\ell}$. 
\end{proof}

\section{Proof of Main Theorem}

In this section, we use our knowledge of patches to prove the main theorem. Let $w$ be a combinatorial pattern for which we have chosen a patching set. It will be useful to introduce a random variable $Y$ to denote the distance between consecutive patches in a randomly chosen sequence.  We will let the the random variable $X$ denote the 
length of the longest $w$-following subsequence of a random sequence from the middle of one patch to the middle of the next that includes the first and last $2m$ elements. 

Ultimately, we will prove Theorem \ref{thm:main} by showing that $L_n^w$ satisfies the central limit theorem of Equation (\ref{e:1}), with $\mu = \mu_X/ \mu_Y$ and $\sigma = \sqrt{\frac{1}{\mu_Y}} \sigma_X$.

Lemma \ref{l:interval} implies that the distribution over the length of the longest subsequence that follows $w$ can be well approximated by the following:  
\begin{enumerate}
\item Begin by waiting for the first occurrence of a patch.  Find the longest subsequence that follows $w$ that ends at the center of the patch and uses the length $2m$ window just before the center of the patch.
\item Repeatedly wait for the next occurrence of a patch and find the longest subsequence between the middle of the previous patch and the middle of the next patch that follows $w$ and uses the entire length $2m$ windows at the beginning and end.  
\item Stop when your original sequence achieves length $n$.  Return the sum of the lengths of subsequences you've found in steps 1 and 2.
\end{enumerate}

In the large $n$ limit neither the contributions to the total length coming from step 1, nor the approximation factor (bounded by the additive amount $n- j_{\ell}$) are large enough to contribute to the limiting distribution. This follows since both quantities are upper bounded by $Y$ which has constant mean and variance, so its contribution when divided by $\sqrt{n}$ limits (as $n\rightarrow \infty$) to  $0$.  Therefore, the driving step for this process is step $2$ which can be seen as repeated sampling from the random variable $X$.  The length is thus the sum of a number of i.i.d. copies of $X$. However the number of copies is not fixed but rather is determined by a partially dependent process, namely the total number of sequence elements drawn.  

Anscombe's Theorem addresses this situation, provided that the distribution of number of copies, $v(n)$ has the property that $\frac{v(n)}{n}$ converges to a constant in probability as $n\rightarrow \infty$ (see for instance \cite[Theorem 2.3]{gut}).  This is indeed the case: the variable $Y$ has finite mean and variance (since the probability that a random sequence of length $\const m$ is an element of $B$ is constant) and this is enough for $\frac{v(n)}{n}$ to converge in probability to $\frac{1}{\mu_{Y}}$ as we now show.  We are interested in showing that for any  constant $\epsilon$ the probability that $|\frac{v(n)}{n} - \frac{1}{\mu_Y}| > \epsilon$ goes to $0$ as $n\rightarrow \infty$.  This condition holds if either the sum of $\lfloor (\frac{1}{\mu_Y} + \epsilon )n \rfloor$ independent copies of $Y$ sum to less than $n$ or that $\lceil (\frac{1}{\mu_Y} - \epsilon )n \rceil$ independent copies of $Y$ sum to more than $n$.  Both events would represent the i.i.d.~sum of copies of $Y$ deviating from their expected value by $O(\epsilon) \sqrt{n}$ standard deviations and normality of the limiting distribution implies their probability for large $n$ is bounded by $e^{-O(\epsilon)^2 n}$ which vanishes as $n \rightarrow \infty$.  With the conditions of Anscombe's Theorem verified and setting $\mu _X$, $\sigma^2_X$ to be the mean and variance of $X$, the application of Anscombe's theorem yields the central limit result stated in Theorem \ref{thm:main} with $\mu = \mu_X/ \mu_Y$ and $\sigma = \sqrt{\frac{1}{\mu_Y}} \sigma_X$.

\section{Algorithms}
We now consider the algorithmic problem of finding long subsequences of a given sequence $s$ that follow a specified combinatorial pattern $w$.
In this section we show two results about the resources needed for this problem.
The first treats the $r=2$ case, that is, up/down patterns.
The second provides an approximate solution in the general case.

In what follows we suppose that the number of bits required to express the elements of a sequence $s$ of length $n$ is $O(\log n)$,
as if $s$ were a permutation of $\{1,\ldots,n\}$ (rather than a point of $[0,1]^n$).

\begin{theorem} \label{l:1}
Given a non-constant up/down pattern $w$ of length $k$ and an input sequence $s$ of length $n$ whose elements are encoded in $O(\log n)$ bits,  there is an $O(n\log n)$ time algorithm that generates a longest $w$-following subsequence of the sequence $s$.  The algorithm uses $O( n \log n)$ bits of memory.  If only the length of a longest subsequence is desired, then only $O(\log n)$ bits of memory are needed. 
\end{theorem}

\begin{theorem} \label{l:2}
For any combinatorial pattern $w$ with combinatorial constant $c$, and for any $\epsilon >0$, there is a $O(n^{\ell c}\log n)$ time algorithm with $\ell= \lceil \frac{1}{\epsilon} \rceil$ that takes an input sequence $s$ of length $n$ (with elements encoded in $O(\log n)$ bits) and generates a $w$-following subsequence of $s$ that has length greater than $(1- \epsilon)$ times the longest $w$-following subsequence of $s$. 
\end{theorem}

Our approach for the first algorithm is in the spirit of a dynamic program.  We step through the given sequence while building a small set of candidate subsequences.  For each $i$, our set of candidate subsequences $S_i$ has the guarantee that it includes a subsequence that matches an optimal subsequence up to position $i$.  As a result, $S_n$ will be guaranteed to contain a longest subsequence that follows the pattern.  We shall make sure our sets $S_i$ do not get too big, thus guaranteeing an efficient algorithm.

\begin{proof}[Proof of Theorem \ref{l:1}]
Starting with $S_0 = \{ \}$, we iterate over $i=1,2, \dots n$ at each time updating our candidate set of subsequences $S_{i-1}$ to $S_i$.  The process of updating has two steps: adding and trimming. Denote the input sequence by $(a_1,\ldots,a_n)$.  For the adding step, we let $S'_i:= S_{i-1} \cup (S_{i-1} \cat (a_i))$ where $S_{i-1}\cat(a_i)$ denotes the set formed from the sequences $s\cat (a_i)$, $s\in S_{i-1}$. Not all the sequences in $S'_i$ will necessarily follow the pattern; in the trimming step we first remove the sequences in $S'_i$ that don't follow the pattern.   We then perform two further trimming steps to end up with the updated set $S_i$:

\begin{itemize}
\item {\it Length trimming.} We remove those subsequences that are more than $c+k$ shorter than the longest remaining element of $S'_i$. 

\item{\it Redundant trimming.}  For each fixed length, we remove all remaining subsequences except those whose final element is maximal or minimal among sequences of that length.  
\end{itemize}

This results in a set containing at most two sequences of each of at most $c+k+1$ possible lengths, hence $|S_i| \le 2(c+k+1)$ independent of $i$.

We claim now that if
$S_{i-1}$ contains an element $t$ that matches an optimal subsequence up to position $i-1$ then $S_i$ contains an element that matches an optimal subsequence up to position $i$.  Note that after the addition step the set $S'_i$ contains both $t$ and $t\cdot(a_i)$, one of which matches an optimal sequence up to position $i$.

The combinatorial property of the up/down pattern ensures that the {\it length trimming} step only removes subsequences that could not match an optimal subsequence up to position $i$.  For if  $t$ and $u$ are subsequences of  $a[1,i]$ with $|t| > |u| +c+k$, and a subsequence $v$ of $a[i+1,n]$ is such that $u v$ follows the pattern, then $v[j, |v|]$ follows $w$ for some $1\leq j\leq k$ and the merge of $t$ and $v[j,|v|]$ will also follow the pattern and be at least as long as $|t| + |v| - j + 1 - c > |u| + |v|$.  This ensures that $u$ cannot match an optimal subsequence up to position $i$.

Now suppose a sequence $u$ of length $l$ is removed in the \emph{redundant trimming} step.  This means there are sequences $u^-,u^+$ in $S_i$ of the same length $l$ and satisfying $u^-_l < u_l < u^+_l$.  Therefore if $u$ matches an optimal sequence up to position $i$ then so does either $u^-$ or $u^+$, depending on whether the entry $W_l$ of the pattern (with the subscript $l$ interpreted mod $k$ as usual) is ``up'' or ``down.''

The claim follows, so $S_n$ contains an optimal subsequence.

The adding and trimming steps in moving from $S_i$ to $S_{i+1}$ take a constant number of operations.  Each of these operations involve comparisons of elements with $O(\log n)$ bits which leads to the $O(n\log n)$  runtime.
\end{proof}

\begin{proof}[Proof of Theorem \ref{l:2}]

Our strategy is as follows.  We start at the beginning of the sequence and find the smallest $i_1$ such that there exists a subsequence $t_1$ of $s[1,i_1]$ of length $\ell c$ ($\ell$ a constant to be specified later) that follows the pattern.  We then repeat this process and find the smallest $i_2$ such that there exists a subsequence $t_2$ of $s(i_1 , i_2]$ of length $\ell c$ that follows the pattern.  We continue this, finding $i_3, \dots, i_b$ and subsequences $t_3, \dots , t_{b}$ until we cannot continue, at which point we have that the maximal subsequence $t_{b+1}$ of $s(i_b, n]$ must be of length some $\ell_0 < \ell c$ which we also compute.  We concatenate these $b+1$ subsequences $t_1, \dots , t_{b+1}$ together and merge them producing a $w$-following sequence $t$ with length at least $|t| \geq b\ell c + \ell_0 - bc$.  

Furthermore, a longest subsequence $a$ of $s$ that follows $w$ cannot have length $L$ bigger than $b\ell c + \ell_0$.  To see this, denote by $\alpha_j$ the index of $s$ that is the location of the $\ell cj$-th element of $a$.  The key observation is that $\alpha_j \geq i_j$ for all $j$ which we now show.  Suppose otherwise and denote by $j_0$, the first $j$ such that $\alpha_{j_0} < i_{j_0}$.  This implies that the sequence $a (\ell c(j_0-1) , \ell c j_0]$ is a length $\ell c$ subsequence that follows $w$ and is contained within $s(\alpha_{j_0-1} , \alpha_{j_0}]$ which would imply that $i_{j_0} \leq \alpha_{j_0}$ which is a contradiction.   We thus establish that $\alpha_b \geq i_b$ which means that the tail of $a$, namely $a(b \ell c, L]$,  must be a  subsequence of $s(i_{b}, n]$ that follows $w$.  This tail therefore has length upper bounded by $\ell_0$ and therefore the length of $a$ must be upper bounded by  $b\ell c +\ell_0$.

Consequently, we have bounded the length $|t|$ of our approximating subsequence $t$  to the interval $[b(\ell-1) c + \ell_0, b\ell c + \ell_0]$.  
Since $\frac{bc }{ b\ell c + \ell_0} < \frac{ 1}{\ell}$, choosing $\ell= \lceil \frac{1}{\epsilon} \rceil$ yields that length of the sequence $t$ lies in the interval $[(1-\epsilon )L, L]$.

To complete the argument, we describe how to find the $i_j$ and $t_j$.  To do this we step along in $i$ and test all the length $\ell c$ subsequences of $s[1,i]$ until we find a $t_j$ that follows the pattern; we set $i_j$ to be the current value $i$.  This is an exhaustive search over at most $n$ choose $\ell c$ subsequences and takes  $O(n^{\ell c})$ comparison steps which translates to $O(n^{\ell c} \log n)$ time.
\end{proof}

\section{A dynamical system}\label{sec:UUD}

In this section, we show that for some patterns $w$, the process of determining the longest $w$-following subsequence can be modeled by a dynamical system on a union of closed connected subsets of $\R^n$. Let $I=[0,1]$ be the unit interval, and $\N$ denote the nonnegative integers.

\begin{example}\label{ex:UD} 
    Consider the pattern $UD$. We codify an online algorithm, different from that in the previous section, for determining the length of the longest $UD$-following subsequence. 
    We first introduce the {\it state space} 
    $$X'=(V_1 \cup V_2 )\times \N,$$
    where $V_1=\{U\}\times I$ and $V_2=\{D\}\times I$.
    A point $(U,x,t)\in X'$ is intended to represent that condition that we are looking for an ``up,'' and we wish this ``up'' to extend a sequence that ends in the number $x$ and whose length is $t$. The point $(D,x,t)$ has a similar interpretation when we are looking for a down. We then have a process
    $$
    \gamma': X'\times I \rightarrow X'
    $$
defined by 
\begin{align*}
    ((U, x, t), s) &\mapsto 
    \begin{cases}
        (D,s, t+1) & s>x \\
        (U,s, t) & s<x.
    \end{cases}\\
    ((D, x, t), s) &\mapsto 
    \begin{cases}
        (D,s, t) & s>x \\
        (U,s, t+1) & s<x.
    \end{cases}\\
\end{align*}
For any sequence $s=(s_1, \dots, s_n)$, we start in the state $(U, s_1, 1)$ and run the process.  One can check that the final value of $t$ is the length of the longest $UD$-following subsequence of $s$.

One can also run this process without the $t$ coordinate of $X'$. That is, let $X=V_1 \cup V_2  $, and consider the process $\gamma:X\times I \rightarrow X$, obtained from $\gamma$ by ignoring the $t$ coordinate. Then $\gamma$ is a process on the compact space $X$ (the disjoint union of two closed intervals), and one can calculate that the unique stationary distribution $p$ on $X$ is given by 
\begin{align*}
    p(U,x) &= 1-x,\\
    p(D,x) &= x.
\end{align*}
From this stationary distribution, when the next symbol $s_*$ is read, the probability of transitioning from $U$ to $D$ or vice versa is $2/3$. Thus if $s$ is a sequence of length $n$, where $n$ is large, then the expected length of the longest $UD$-following subsequence is approximately $(2/3)n$, in agreement with Stanley \cite{stanley-mich}. 
\end{example}

\begin{example}\label{ex:UUD} 
    Consider the pattern $UUD$. Now the algorithm is a bit more complicated; it turns out that when we are looking for the second $U$, which we denote as $U_2$, we must remember two previous values $(x,y)$, where $x<y$.  For this purpose, we introduce the triangle $T=\{(x, y)\in I\times I | x<y\}$. 
    This time the state space $X$, without a counting parameter $t$, is the union of two segments and a triangle:
    $$
    X=V_1\cup V_2\cup V_3,
    $$
    where $V_1=\{U_1\}\times I$, $V_2=\{U_2\}\times T$, $V_3=\{D\}\times I$. The process
    $\gamma: X\times I \rightarrow X$ is given by
\begin{align*}
    ((U_1, x),s) &\mapsto 
    \begin{cases}
        (U_2,(x,s)) & s>x \\
        (U_1,s) & s<x.
    \end{cases}\\
    ((U_2, (x,y)), s) &\mapsto 
    \begin{cases}
        (D,s) & s>y \\
        (U_2,(x,s))& x<s<y\\
        (U_2,(s,y)) & s<x\\   
    \end{cases}\\
     ((D, x), s) &\mapsto 
    \begin{cases}
        (D,s) & s>x \\
        (U_1,s) & s<x.
    \end{cases}\\
\end{align*}
For any sequence $s=(s_1, \dots, s_n)$, we start in the state $(U_1, s_1)$ and perform this process. The length of the longest $UUD$-following subsequence can be found by counting the number of times that the first coordinate changes, i.e., the number of times the system changes connected components.  Figure \ref{fig:UUD} shows the component-level structure of this dynamical system.

\begin{figure}[b]
    \centering
    \begin{tikzpicture}[scale=1]
        \draw[very thick] 
            (3,0) -- (3,1);
        \draw[very thick]
            (0,3) -- (1,3);
        \draw[fill] (1,0) -- (1,1) -- (0,1) -- cycle;

        \draw (3,-.5) node {$D$};
        \draw (-.5,3) node {$U_1$};
        \draw (-.5,1) node {$U_2$};

        \draw [->] (.5,2.5) -- 
            node[anchor=east]{$s>x$}         
            (.5,1.5);
        \draw [->] (1.5,.5) -- 
            node[anchor=north]{$s>y$}         
            (2.5,.5);
        \draw [->] (2.5,1.5) -- 
            node[anchor=south west]{$s<x$} 
            (1.5,2.5);

        \draw[->] (3.5,.25) arc (-135:135:.6 and .3);
        \draw[->] (.75,3.5) arc (-45:225:.3 and .6);
        \draw[->,rotate=45] (.25,.4) arc (45:315:.5 and .1);
        \draw[->,rotate=45] (.25,-.2) arc (45:315:.5 and .1);

        \draw (5,0) node {$s>x$};      
        \draw (1.5,4) node {$s<x$};      
        \draw (-1.25,0) node {$s<x$};     
        \draw (-1.2,-.6) node {$x<s<y$};     

    \end{tikzpicture}
    \caption{The dynamical system $X$ for the pattern $UUD$.}
    \label{fig:UUD}
\end{figure}

As in Example \ref{ex:UD}, there is a unique stationary distribution on $X$ for the process $\gamma$.  This distribution solves an explicit system of integral equations, but is largely a mystery to us. Numerically, we have found that from this distribution the chance of transitioning is approximately $\delta=0.577447517$.   Thus for large $n$, the expected length of the longest $UUD$-following subsequence of $s$ is approximately $\delta n$. We calculated this numerical approximation for $\delta$ by direct simulation.  More efficiently, one can iterate the process analytically, starting with a point mass, keeping track of the resulting distribution on $X$, which is given by polynomial functions on the $V_i$.  Convergence appears to be rapid; for example, after evolving for 70 steps, we found that consecutive distributions have an $L^2$ distance of approximately $10^{-15}$. 
\end{example}

It may be interesting to study families of patterns as well.  For instance denote by $\delta_k$ the mean associated to the pattern $U^kD$, so $\delta_1=2/3$ and $\delta_2$ equals the value of $\delta$ from Example 2.  Is $\delta_k$ decreasing in $k$?  Does $\delta_k \to 0$?

\begin{example}\label{ex:UUDD}
Consider the pattern $UUDD$. For this pattern, there is an algorithm, similar to the $UUD$ algorithm, whose state space $X$ is the union of two intervals and two triangles. When looking for the first of the two ups or the first for the two downs, one keeps track of a single previous value; when looking for the second of the two ups or the second of the two downs, one keeps track of two previous values. This process has a form similar to that of Example \ref{ex:UUD}, and we omit a detailed description. Again, the stable distribution, which we approximated numerically, is a mystery. We can report an approximate transition probability of $\mu \approx 0.561$. Hence for large $n$, in a sequence of length $n$ the longest $UUDD$-following subsequence has expected length approximately $0.561n$.
\end{example}

\begin{question}
    Is every nonconstant up/down pattern modeled by such a dynamical system on a union of simplices?  Are there dynamical systems modeling combinatorial $r$-patterns for $r>2$?
\end{question}

\section{Further questions} \label{sec:questions}

We conclude with a discussion of some additional open questions.
In light of Theorem \ref{thm:main}, the following updates Stanley's question.
\begin{question}
    Are there patterns $w$ such that the expected length of the longest $w$-following subsequence of a random length $n$ sequence is neither $\Theta(n^{1/2})$ nor $\Theta(n)$?
\end{question}
It seems possible that the pattern $w_1$ defined in \eqref{eq:squeeze} below has length $o(n^{1/2})$.
However, we have the following general lower bound, which applies to any pattern $w$, whether combinatorial or not.

\begin{theorem}[General lower bound] \label{thm:genlower}
    For any pattern $w$, and any fixed $\epsilon>0$, we have $\mu^w_n \ge n^{1/2-\epsilon}$.
\end{theorem}

\begin{proof} We view the sequence $s=(s_1,\ldots,s_n)$ of length $n$ as uniformly random in $[0,1]^n$.
Set $m=n^{1/2 -\epsilon}$ and let $\sigma\in S_m$ be a permutation that follows $w$.
We argue that there is high probability that there exists a length $m$ subsequence of $s$ with order $\sigma$.

Break the interval $[0,1]$ into the $m$ ``bins'' $[0,1/m), [1/m,2/m), \cdots,[1-1/m,1]$.
We refer to $[(i-1)/m,i/m)$ as bin $i$.

Let $N=n/m = n^{1/2+\epsilon}$.  The chance that (at least) one of the first $N$ elements of $s$ lies in bin $\sigma(1)$ is $1-(1-1/m)^N$.
This is also the probability of each (similar but) independent event $E_j$ (for $1\le j \le m$) that the element $s_i$ lies in bin $\sigma(j)$, for at least one $i$ in the interval $((j-1)N , jN]$.

Note that if all events $E_j$ occur then $s$ contains a $w$-following subsequence of length $m$.

The probability that all $E_j$ occur is  
$$\left(1-\left(1-1/m\right)^N\right)^m \geq 1 - m (1-1/m)^N
= 1-m (1-1/m)^{m\cdot n^{2\epsilon}}$$
which is asymptotically $1-me^{-n^{2\epsilon}}$ which, as $\epsilon$ is fixed, tends to 1.

Thus with high probability the length of the longest $w$-following subsequence of $s$ is at least $m$, and the result follows.
\end{proof}

In fact the previous argument shows slightly more, namely that $\mu^w_n \ge f(n)$ if $f(n)$ is any function satisfying 
$\log f - n/f^2 \to -\infty$, for example $f(n)=n^{\frac 12 - \frac{\log\log n}{\log n}}$.

\subsection{Drift}
We preface another collection of questions by mentioning that this work grew in part out of our discovery of a phenomenon called \emph{drift} in a previous project \cite{tc09}.
A (perhaps overly) formal definition appears there but for the present purpose, a pattern $w$ of length $k$ has ``upward drift'' in coordinate $i\in \Z_k$
if every $w$-following sequence $s=(s_1,\ldots)$ satisfies $s_j < s_{j+k}$ for all $j\equiv i \mod k$.  Examples include the increasing pattern $U$ and the 3-pattern
$W_1=\{(1\ 3\ 2)\}, W_2=\{(2\ 1\ 3)\}$ of length 2, which has upward drift in both coordinates.  Downward drift is defined similarly.  A pattern can drift  differently in different coordinates, 
e.g. 
\begin{equation}\label{eq:squeeze}
w_1 \mbox{ defined by } W_1=\{(1\ 3\ 2)\}, W_2=\{(3\ 1\ 2)\},
\end{equation}
and a pattern can drift in some coordinates and not others, e.g. the length 4 pattern with 
$W_1=\{(1\ 3\ 2 )\}, W_2=\{(2\ 1\ 3)\}, W_3=\{(2\ 3\ 1)\}, W_4=\{(2\ 1\ 3)\}$.
A pattern is called \emph{driftless} if there is no drift in any coordinate; one example is the alternating $UD$ and another is the 3-pattern 
\begin{equation}\label{eq:dnotc}
w_2 \mbox{ defined by } W_1=\{(1\ 2\ 3)\}, W_2=\{(2\ 3\ 1)\}, W_3=\{(3\ 1\ 2)\}
\end{equation}
of length 3.

Our primary interest in this article is the distribution of the lengths of longest $w$-following subsequences.
Patterns with drift have a mean length at most $O(\sqrt{n})$ by comparison with the increasing case, 
and our initial suspicion was that driftless patterns might resemble the alternating case, i.e.~display linear growth in $n$.

Consider the following statements about a pattern $w$:

\begin{description}
\item[(C)] $w$ is combinatorial;
\item[(G)] the distribution of $1/n$ times the length of the longest $w$-following subsequence of a random sequence of length $n$ tends to a Gaussian as $n$ grows;
\item[(L)] the expected length of the longest $w$-following subsequence of a random sequence of length $n$ is linear, i.e.~is asymptotic to 
$\mu n$ for some constant $\mu>0$, as $n$ grows;
\item[(D)] $w$ is driftless.
\end{description}

These conditions are successively weaker:
\begin{equation}\label{implies}
\mbox{\C} \implies\mbox{\G} \implies\mbox{\L} \implies\mbox{\D}
\end{equation}
The first implication is our main theorem, the second is immediate, and the third follows from the Baik-Deift-Johansson theorem.

For 2-patterns, which we have also called ``up/down'' patterns, it is also true that \D\ implies \C.
Thus for 2-patterns of any length these categories are all equivalent.  More specifically, constant 2-patterns (such as the increasing $U$ and the decreasing $D$) fall into none of the types \C, \G, \L, \D, and all other 2-patterns are of type \C, \G, \L, \D.
In particular, the dichotomy we imagined is a reality for up/down patterns.

If $r\ge3$ then \D\ does not imply \C, as we now show.  Suppose that $w$ is combinatorial and $s$ and $t$ are $w$-following sequences such that all the elements $s$ are greater than all the elements of $t$. Then the combinatorial condition applied to $st$ implies that $w$ must satisfy a condition that we call {\em patching down}.
An $r$-pattern $w$ of length $k$ is said to patch down if there exists $j\in \Z/k\Z$ such that for each $1\le i \le r$ there is a permutation in $W_{j+i}$ whose last $i$ entries are (in some order) the numbers $1,2,\ldots,i$.  Likewise the combinatorial condition applied to $ts$ leads to a {\em patching up} condition on  $w$:  an $r$-pattern $w$ of length $k$ is said to patch up  if there exists $j\in \Z/k\Z$ such that for each $1\le i \le r$ there is a permutation in $W_{j+i}$ whose last $i$ entries are (in some order) the numbers  $r,r-1,\ldots,r-i+1$.
 The driftless pattern $w_2$ defined in \eqref{eq:dnotc}
patches down but not up, hence cannot be  combinatorial.

It follows that at least one of the implications \eqref{implies} is not reversible.

\begin{question} 
    Are any of the above implications reversible?
\end{question}

We suspect but have not proved that \D\ does not imply \L.
The suspicion comes from the pattern $w_2$ defined in \eqref{eq:dnotc} above.  Because it doesn't patch up, one could have a long subsequence following $w$ followed by a longer and higher subsequence following $w$ and so on.  This behavior is similar to the decreasing pattern $D$ which fails \L.

In general \L\ also does not imply \G, though the implication may hold for simple patterns.  
Consider the following 4-pattern of length 1:  
\begin{equation*}
w_3=\left\{\{ (1234), (2341), (3412), (4123) , (4321), (3214), (2143), (1432) \}\right\}.
\end{equation*}  
A subsequence following $w_3$ either follows the cyclic permutations of $(1234)$ or else follows the cyclic permutations of $(4321)$ but cannot jump between the two.  Each of these patterns separately is combinatorial, so the length of the longest $w_3$-following subsequence is distributed as the maximum of two (correlated) Gaussians, which is not Gaussian.
Again, however, it is possible that \L\ implies \G\ for simple patterns.

The claim in the preceding paragraph that various patterns are combinatorial is somewhat tedious to verify.
If one knows $c$, then something like this can always be done by exhaustive case analysis; in particular there is an algorithm to find the smallest possible $c$ if one knows that $w$ is combinatorial.  We mentioned earlier that we know of no example where $c$ needs to be larger than $\max\{k,r\}$, but we have not proved this.

\begin{question}
    How hard is it to decide whether a given pattern is combinatorial?
\end{question}

Finally, we do not know whether \G\ implies \C.

\subsection{Computational issues}
Algorithmically, the time and space complexity of the problem of finding the longest $w$-following subsequence (or its length, or approximate versions) seem to be related to the properties isolated above.  For example, Theorem \ref{l:2} says that if $w$ is combinatorial, then there is a log-space algorithm for approximating the length of the longest $w$-following sequence. In fact, the algorithm given there is a {\it one-pass} algorithm, meaning each bit of the input is read exactly once in order from left to right. The situation for the pattern $U$ is quite different. Finding the length of the longest increasing subsequence, with the restriction of only reading the sequence once,  requires one to keep track of a linear, rather than constant, number of potential candidate sequences \cite{streaming-lis}. For patterns $w$ with drift, one would thus expect there does not exist a log-space one-pass algorithm for finding the length of the longest $w$-following subsequence.

We introduce the condition
\begin{description}
\item[(A)] There exists a one-pass log-space algorithm for finding the longest $w$-following subsequence.
\end{description}
One can then ask how \A\ fits into the {\bf CGDL} hierarchy.  Variants of \A\ might also be interesting. For example, one could drop the one-pass restriction, requiring only that the algorithm be log-space. Or, in the spirit of Theorem \ref{l:2}, one could ask for approximate algorithms instead of exact ones.

\subsection{Circles}
We consider the related setup in which a sequence of points is given on the unit circle, rather than the unit interval, and we seek the longest subsequence following a given \emph{circular pattern}, which is defined just like a pattern except we
look at the \emph{circular order} of each window of size $r$ from the sequence.  The circular order of $r$ points on the circle is the order in which the points are encountered as we travel counterclockwise on the circle, starting at the first point. For instance, for a circular sequence $s$ to follow the circular 3-pattern $\{(1 2 3)\}$ we require that for each $i$, if we make one counterclockwise revolution around the circle starting at $s_i$, we encounter $s_{i+1}$ before $s_{i+2}$. 

The circular order of any sequence on a circle will have first entry 1, and on a circle every sequence follows the 2-pattern $(12)$, just as on an interval every sequence follows the 1-pattern $(1)$.

What we noticed is that the length associated to the circular 3-pattern $\{(1 2 3)\}$ has mean $\mu=2/3$ and variance $\sigma^2=8/45$, the same values arising from alternating sequences in an interval.  A direct proof that the two should agree would be interesting, and may lead to an understanding of additional circular patterns.

\section{Appendix:  Construction of $\tau$}

Our goal is to prove Lemma \ref{l:tau} stated below.  We'll make use of the following lemma:
\begin{lemma} [Extending a tracking sequence] \label{l:addon} 
If a sequence $t$ of length $2m$  tracks $\sigma \in S_m$, then there exists a sequence $u$ of length $m$ such that $tu$ also tracks $\sigma$.
\end{lemma}
\begin{proof}
Let $\tilde{t}$ be a tie break for $t$ and let $\alpha$ be the order of $\tilde{t}$ and note that $\alpha$ tracks $\sigma$.  Since  $t(m, 2m]$ also tracks $\sigma$ we can create $u$ so that $t(m, 2m] \cat u$ tracks $\alpha$ by defining any strictly increasing continuous map $f: [1, 2m] \rightarrow \mathbb{R}$ with the property that $f(\alpha_i) = t_{m+i}$ (whichh can be done because the order of $\alpha[1,m]$ and $t(m, |t|]$ are both $\sigma$) and then setting $u[i]= f(\alpha_{m+i})$. 

To show such a $tu$ tracks $\sigma$ we must show  the order of $(t\cat u)[i, i+m )$ is $\cyc^{i-1}(\sigma)$ for $1 \leq i \leq 2m+1$.  This holds for $i\leq m$ because $t$ tracks $\sigma$.   For $i \geq m+1$ this follows from the fact that the order $(t\cat u)[i,i+m)$ is the order of $\alpha[i-m, i)$ which by virtue of $\alpha$ tracking $\sigma$ is the order of $\cyc^{i-m-1}(\sigma)= \cyc^{i-1}(\sigma)$.
\end{proof}

Fix a combinatorial $r$-pattern $w$ of length $k$ and let $c$ be the corresponding combinatorial constant.  Consider the directed graph $H=H_w$ with vertices $V=\{ (\sigma, i): \sigma \in W_i,\ i\in\Z_k \}$ and directed edges connecting $(a,i)$ to $(b, i+1)$ whenever the order of the last $r-1$ elements of $a$ equals the order of the first $r-1$ elements of $b$.  Given a vertex $v \in V$ define the {\it loop orbit of $v$} to be the set of vertices of $H$ such that there is a loop\footnote{By a \emph{loop} we mean a sequence of vertices beginning and ending with the same vertex such that there is a directed edge from each vertex to the next.  This is sometimes called a \emph{circuit} in graph theory literature.} in $H$ containing both $v$ and $w$. This is also the set of vertices of the strongly connected component of $H$ that contains $v$. 

Note that if the pattern $w$ is simple, then $H$ is a directed cycle of length $k$.  Some of the arguments that follow are more straightforward in this case.

For $w$-following sequences $s$ and $t$ we'll say that {\it $s$ can reach $t$} if there is a sequence $u$ such that $sut$ is $w$-following and $|su|$ is a multiple of $k$.

\begin{lemma}[Construction of $\tau$] \label{l:tau} Given a combinatorial $r$-pattern $w$ of length $k$ with combinatorial constant $c$ there exists an $m\geq \max (4c,4k,r)$ that is a multiple of $k$ and a permutation $\tau \in S_m$ with the following properties:
\begin{enumerate}
\item  Any sequence $s$ that tracks $\tau$ also follows $w$.
\item When $\tau$ is interpreted as a path on $H$, it passes through all the vertices in the loop orbit of $(\sigma', 1)$ where $\sigma'$ is the order of $\tau[1,r]$.
\item If there exists a $w$-following sequence $u$ for which the first $2m$ and last $2m$ elements of $u$ both track $\tau$, then when $u$ is interpreted as a path on $H$ it only passes through vertices in the loop orbit of $(\sigma',1)$.
\end{enumerate}
\end{lemma}

\begin{proof} Among vertices of $H$ with second coordinate $1$, let $(\sigma,1)$ be a vertex with a maximal size loop orbit.  Let $L$ be a loop in $H$ starting at $(\sigma,1)$ that contains all of $(\sigma, 1)$'s loop orbit and is at least of length $\max (2c,2k,r)$.  Create a permutation $\eta$ such that the order of $\eta[i, i+r)$ corresponds to the first coordinate of the $i$th vertex in $L^{4}$ ($4$ consecutive loops around $L$).  By construction, $\eta$ follows $w$ and consequently the combinatorial property says we can find a long subsequence $\zeta$ of the sequence  $\eta^{\ell}$ that also follows $w$, where we will specify the size of $\ell$ in a moment.  Our goal is to use the pigeon hole principle to show there is a sufficiently long window of $\zeta$ of the form $ABA$ where the two $A's$ represent identical length $r$ strings and the length of $AB$ is $dk$ for some integer $d \ge 2$.    Let $S$ be the set of length $r$ windows of $\eta$ ($A$ will turn out to be some element of $S$), i.e.  
\[S=\{\eta[i,i+r): \ 1\leq i \leq |\eta| - r +1\}. \]
Looking at windows of size $r$ in $\zeta$, we are guaranteed to see many windows that lie in $S$ since the combinatorial property used to form $\zeta$ keeps all but possibly the first and last $c$ elements of each $\eta$ untouched.  Of the original $|\eta| -r +1$ windows of length $r$ in each $\eta$, at least $|\eta| -r +1 - 2c$ of them remain within $\zeta$.  Combining this lower bound across each of the $\ell$ copies, we can guarantee that we see at least $\ell(|\eta| - r+1- 2c)\geq \ell (4k)$ elements of $S$ when looking at windows of length $r$ of $\zeta$. Choosing $\ell= |S|$ so that $\ell(|\eta|- r+1 -2k) \geq 4k|S| \geq k|S| +1$ implies that we can find two locations $i_0, j_0$, $1\leq i_0 \leq j_0 \leq |\zeta|-r+1$ such that
\begin{itemize}
    \item $\zeta[i_0, i_0+r)=\zeta[j_0, j_0+r) \in S$
    \item $j_0 -i_0= 0 \mod k$
    \item $j_0-i_0 \geq |\eta|-2c .$
\end{itemize}
The key observation is that the section of $\zeta$ between $i_0$ and $j_0 + r-1$ represents a loop in $H$ since $\zeta$ follows $w$ and the starting and ending vertex are the same as guaranteed by the first two bullet points above. It follows that repeating $\zeta[i_0, j_0)$ will trace around this loop over and over.  We define $\rho$ to be the length $m=j_0-i_0$  order corresponding to a tie break of a shifted version of $\zeta[i_0, j_0)$ so that as a loop on $H$ the first element has second coordinate $1$: 
\[\rho = \mbox{the order of a tie break of } \cyc^{\star}\left(\zeta[i_0,j_0)\right), \]
where $\star=k\lfloor\frac {j_0} {k}\rfloor$ is the largest multiple of $k$ less than or equal to $j_0$.   We conclude that any power of $\rho$ follows $w$ and thus any sequence that tracks $\rho$ will follow $w$.   

We now show that taking $\tau$ to be a well chosen shift (possibly no shift at all) of $\rho$ will yield the results of the lemma.  Define $I$ to be the set of indices $i$ such that $\cyc^i(\rho)$ follows $w$; as noted, our construction ensures that $0\in I$. 
We'll say $(\cyc^i(\rho))^2$ reaches  $(\cyc^j(\rho))^2$ if there is a sequence $s$ that tracks $(\cyc^i(\rho))^2$ and a sequence $t$ that tracks $(\cyc^j(\rho))^2$ so that $s$ reaches $t$.  Set $\tau= \cyc^{i}(\rho)$ where $i \in I$ is chosen so that $(\cyc^{i}(\rho))^2$ reaches the minimal number of  $(\cyc^j(\rho))^2$.   Since the reaching property is transitive, this minimality ensures that if $\tau^2$ can reach $(\cyc^j(\rho))^2$ then $(\cyc^j (\rho))^2$ can also reach $\tau^2$. 

By construction our choice of $\tau$ satisfies statement (1) of the lemma.
To prove statement (2) we rely on the fact that our choice of $\tau$ was long.  Recall that $|\tau|=m=j_0-i_0$ was at least $|\eta|-2c$ (third bullet point above).  Since the middle two $L$ loops of each copy of $\eta$ that are used to form $\zeta$ are untouched, the length of $\tau$ guarantees that it contains a window that originally corresponded traveling at least once around $L$.  However, because $\tau$ may be shifted, when tracking the corresponding path of $\tau$ in $H$ we cannot conclude the path includes $L$ but rather only a loop with the same number of distinct vertices as $L$.  When we combine this with the fact that $(\sigma,1)$ was chosen so that its loop orbit was of maximal size, we conclude that the corresponding path of $\tau$ in $H$ must include all elements of the loop orbit of $(\sigma',1)$ where $\sigma'$ is the order of $\tau[1,r]$. This establishes the second condition in the lemma. 

For (3), we begin with a $u$ that is $w$-following and such that $u[1,2m]$ and $u(|u| - 2m , |u|]$ track $\tau$.  We now show that we can construct an extension of $u$ that will provide a witness to the fact that $\tau^2$ reaches $\cyc^{|u|}(\tau)^2$. By Lemma \ref{l:addon} we can find a sequence $v$ of length $m$ so that the length $3m$ sequence $u(|u|-2m, |u|]v$ tracks  $\tau$.  This implies that the extended sequence $uv$ also follows $w$: the only condition that needs to be verified is that the order of $uv[i, i+r)$ lies in $W_{i}$ for $i\geq |u|-r+1$.  This holds since for such $i$ the order of $uv[i, i+r)$ is the order of $uv[i-m, i-m +r)$ (because the end of $uv$ tracks $\tau$) which is already known to be in $W_{i-m}= W_{i}$.
(Recall that subscripts of $W$ are interpreted mod $k$.)
Let $p = |u| \mod k$ and define the truncated sequence $u'= u v[1, m - p]$ of length a multiple of $k$ that also follows $w$.  Since the beginning $2m$ elements track $\tau$ and the final $2m$ elements track $\cyc^{m-p} (\tau)$ we have that $\tau^2$ reaches $(\cyc^{m-p} (\tau))^2$.  Our choice of $\tau$ then ensures that $(\cyc^{m-p}(\tau))^2$ also reaches $\tau^2$ and consequently establishes that the path in $H$ corresponding to $u'$ only includes loop vertices of $(\sigma', 1)$ where $\sigma'$ is the order of $\tau[1,r]$. This establishes the last condition of the lemma.
\end{proof}

\end{document}